\documentclass[12pt,reqno]{amsart}

\usepackage{amssymb}
\usepackage{tikz}
\usetikzlibrary{arrows}

\numberwithin{equation}{section}

\def\ZZ{{\mathbb Z}}

\def\NN{{\mathbb N}}
\def\FF{{\mathbb F}}

\newtheorem{lemma}{Lemma}[section]
\newtheorem{corollary}[lemma]{Corollary}
\newtheorem{theorem}[lemma]{Theorem}
\newtheorem{proposition}[lemma]{Proposition}

\theoremstyle{definition}
\newtheorem{definition}[lemma]{Definition}
\newtheorem{remark}[lemma]{Remark}
\newtheorem{example}[lemma]{Example}

\title[Lattice paths with given number of turns and $\Gamma$-semimodules]{Lattice paths with given number of turns and semimodules over numerical semigroups}

\author{Julio Jos\'e Moyano-Fern\'andez \and Jan Uliczka}
\address{Institut f\"ur Mathematik, Universit\"at Osnabr\"uck. Albrechtstr. 28a, \newline D-49076 Osnabr\"uck, Germany. \emph{Email adresses}: {\tt jmoyanof@uos.de} {\rm and}  {\tt juliczka@uos.de} }

\begin{document}


\maketitle

\begin{abstract}
Let $\Gamma=\langle \alpha, \beta \rangle$ be a numerical semigroup. 
In this article we consider several relations between the so-called $\Gamma$-semimodules and lattice paths from $(0,\alpha)$ to $(\beta,0)$: we investigate  
isomorphism classes of $\Gamma$-semimodules as well as certain subsets of the set of gaps of $\Gamma$, and finally syzygies of $\Gamma$-semimodules. In particular we compute the number of $\Gamma$-semimodules which are isomorphic with their $k$-th syzygy for some $k$.
\end{abstract}

\section{Introduction} \label{section1}

In our paper \cite{mu} we considered Hilbert series of graded modules over the polynomial ring $R=\FF[X,Y]$ with $\deg(X)$ and $\deg(Y)$ being coprime.
The central result was an arithmetic criterion for such a series to be the Hilbert series of some finitely generated $R$--module of positive depth. 
This criterion is formulated in terms of the numerical semigroup generated by $\deg(X)$ and $\deg(Y)$. For reader's convenience we recall some basic vocabulary of this theory here.
\medskip

Let $\Gamma$ be a sub--semigroup of $\NN$ such that the greatest common divisor of all its elements is equal to $1$. Then the set 
$\NN \setminus \Gamma$ has only finitely many  elements, which are called the \emph{gaps} of $\Gamma$. 
Such a semigroup is said to be \emph{numerical}.  
The crucial notion in \cite{mu} was that of a \emph{fundamental couple}: Let $\alpha, \beta >0$ be coprime integers and let $G$ denote the set of gaps of $\langle \alpha,\beta \rangle$. An \emph{$(\alpha,\beta)$--fundamental couple} $[I,J]$ consists of two 
integer sequences $I=(i_k)_{k=0}^m$ and $J=(j_k)_{k=0}^m$, such that 
\begin{enumerate}
\item[(0)] $i_0=0$.
\item[(1)] $i_1, \ldots , i_m, j_1, \ldots , j_{m-1} \in G$ and $j_0,j_m \leq \alpha \beta$.
\item[]
\item[(2)] $\begin{array}{lllll} i_k \equiv j_k &\!\!\!\mod \alpha &~\mbox{and}~& i_k < j_k & ~~\mbox{for} ~~ k= 0, \ldots ,m;\\
        j_k \equiv i_{k+1} &\!\!\!\mod \beta &~\mbox{and}~& j_k > i_{k+1} & ~~\mbox{for} ~~ k= 0, \ldots ,m-1;\\
        j_m \equiv i_{0}  &\!\!\!\mod \beta &~\mbox{and}~& j_m \geq i_{0}. &
       \end{array}$
\item[]
\item[(3)] $|i_k - i_{\ell}| \in G ~~\mbox{for} ~~1 \leq k < \ell \leq m$.
\end{enumerate}

One of the problems considered in this article will be the counting of sets of integers like those appearing in the first position of a fundamental couple. We coin a name for these sets:

\begin{definition}
Let $\Gamma$ be a numerical semigroup.
A set $\{x_0=0, x_1, \ldots, x_n\} \subseteq \NN$ is called  \emph{$\Gamma$-lean} if $|x_i-x_j| \notin \Gamma$ for $0 \leq i <j \leq n$.
\end{definition}

The next two sections deal with objects related to $\Gamma$-lean sets: we begin with isomorphism classes of $\Gamma$-semimodules. In section 3 we consider certain lattice paths; from this section on, $\Gamma$ is restricted to be generated by two elements. In the last sections we turn our attention to the second position of a fundamental couple: we identify the sequences $J$ appearing there with so-called syzygies of $\Gamma$-semimodules and investigate their relation with lattice paths. The process of taking syzygies can be iterated --- our final result allows to compute the number of semimodules $\Delta$ whose $k$-th syzygy is isomorphic with $\Delta$ for some $k$.

\section{Generators of $\Gamma$-semimodules} \label{section2}

Let $\Gamma$ be a numerical semigroup. A $\Gamma$-semimodule $\Delta$ is a non-empty subset of $\NN$ such that $\Delta + \Gamma \subseteq \Delta$.
A system of generators of $\Delta$ is a subset $\mathcal{E}$ of $\Delta$ with
\[
\bigcup_{x \in \mathcal{E}} (x+\Gamma) = \Delta. 
\]
It is called minimal if no proper subset of $\mathcal{E}$ generates $\Delta$.
Note that, since $\Delta \setminus \Gamma$ is finite, every $\Gamma$-semimodule is finitely generated.  

\begin{lemma}
Every $\Gamma$-semimodule $\Delta$ has a unique minimal system of generators.
\end{lemma}

\proof Inductively we construct a sequence $(x_i)$ of elements of $\Delta$ starting with $x_1=\min \Delta$ such that any system of generators has to contain this sequence. If $x_1, \ldots , x_n$ are already constructed but do not generate $\Delta$ we set $x_{n+1}=\min \Delta \setminus \cup_{i=1}^{n} (x_i+\Gamma)$. After finitely many steps we arrive at a system of generators $x_1, \ldots, x_r$, and by construction it is clear that this system is minimal and that any system of generators must contain $x_1, \ldots, x_r$.
\qed

\begin{lemma}
Let $x_1, \ldots, x_r$ be the minimal system of generators of a $\Gamma$-semimodule. Then $|x_i-x_j|$ is a gap of $\Gamma$ for all $i \neq j$. Conversely, any subset $\{x_1, \ldots, x_r\}$ of $\NN$ with this property minimally generates a $\Gamma$-semimodule.
\end{lemma}
\proof We may assume $x_j >x_i$ for $j>i$. Then, by minimality, $x_j - x_i \notin \Gamma$ for all $j>i$. The second assertion is clear since $x_i \notin \cup_{j \neq i} (x_j + \Gamma)$. 
\qed
\medskip

Two $\Gamma$-semimodules $\Delta, \Delta^{\prime}$ are called isomorphic if there is an integer $n$ such that $x \mapsto x+n$ is a bijection from $\Delta$ to $\Delta^{\prime}$. 
For every $\Gamma$-semimodule $\Delta$ there is a unique semimodule $\Delta^{\prime} \cong \Delta$ containing $0$; such a $\Gamma$-semimodule is called normalized.
\medskip

The semimodule $\Delta^{\circ}:=\{x-\min \Delta \mid x \in \Delta \}$ is called the normalization of $\Delta$. It is the unique semimodule isomorphic to $\Delta$ and containing $0$.

\begin{corollary} \label{2C2}
The minimal system of generators of a normalized $\Gamma$-semimodule is $\Gamma$-lean, and conversely, every  $\Gamma$-lean set of $\NN$ minimally generates a normalized $\Gamma$-semimodule. Hence there is a bijection between the set of isomorphism classes of $\Gamma$-semimodules and the set of $\Gamma$-lean sets of $\NN$.
\end{corollary}

\section{Lattice paths and $\langle \alpha, \beta \rangle$-lean sets} \label{section3}

From now on we only consider numerical semigroups with two generators $\alpha < \beta$. In this case there is a connection between $\Gamma$-lean sets and certain lattice paths which allows to deduce a formula for the number of $\langle \alpha, \beta \rangle$-lean sets.

\begin{lemma}[\cite{rosales}, Lemma 1, resp. \cite{mu}, Corollary 3.5] \label{3L1}
(1) Let $e \in \ZZ$. Then $e \notin  \langle \alpha, \beta \rangle$ if and only if there exist 
$k, \ell \in \NN_{>0}$ such that $e=\alpha \beta - k \alpha - \ell \beta$.\\
(2) Any integer $n>0$ has a unique presentation $n=p\alpha\beta-a\alpha-b\beta$ with integers $p>0, 0\leq a<\beta$ and $0 \leq b <\alpha$.
\end{lemma}

This result yields a map $G \to \NN^2$, $\alpha\beta -a\alpha-b\beta \mapsto (a,b)$ which identifies a gap with a lattice point. Since $\alpha \beta -a \alpha -b\beta>0$ the point lies inside the triangle with corners $(0,0), (\beta,0),(0,\alpha)$.

\begin{lemma}[\cite{mu}, Lemma 3.19] \label{3L2}
Let $i_1=\alpha \beta -a_1 \alpha - b_1 \beta, i_2=\alpha \beta -a_2 \alpha - b_2 \beta$ be gaps of $\langle \alpha, \beta \rangle$. Then
the difference $|i_1 - i_2|$ is a gap if and only if $(a_2-a_1)(b_2-b_1) <0$.
\end{lemma}

\begin{corollary}\label{2C1}
Let $\mathcal{E}:=\{0,x_1, \ldots , x_m\} \subseteq \NN$ with gaps $x_i=\alpha \beta -a_i \alpha -b_i\beta$ of $\langle \alpha, \beta \rangle$, $i=1, \ldots, m$ such that $a_1<a_2 < \cdots < a_m$, then $\mathcal{E}$ is $\langle \alpha, \beta \rangle$-lean if and only if  $b_1>b_2 > \cdots > b_m$.
\end{corollary}

Therefore an $\langle \alpha, \beta \rangle$-lean set yields a lattice path with steps downwards and to the right from $(0,\alpha)$ to $(\beta,0)$ not crossing the diagonal, where the points identified with the gaps mark the turns from $x$-direction to $y$-direction. In the sequel those turns will be called ES-turns for short. 

\begin{center}
  \begin{tikzpicture}[scale=0.8]
    \draw[] (0,0) grid [step=1cm](7,5);
    \draw[] (0,5) -- (7,0);
    \draw[ultra thick] (0,5) -- (0,3) -- (1,3) -- (1,2) -- (3,2) -- (3,1) -- (4,1) -- (4,0) -- (7,0);
    \draw[fill=white] (0,5) circle [radius=0.1]; 
    \draw[fill] (1,3) circle [radius=0.1]; 
     \draw[fill] (3,2) circle [radius=0.1]; 
     \draw[fill] (4,1) circle [radius=0.1]; 
      \draw[fill=white] (7,0) circle [radius=0.1]; 
   \end{tikzpicture}
\end{center}
\begin{center}
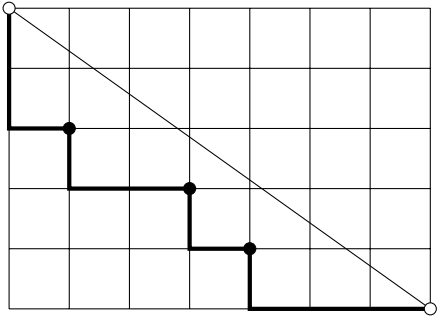
{\small Lattice path for the $\langle 5,7 \rangle$-lean set $\{0,9,6,8\}$.}
\end{center}
\medskip

(Note that in the first component $I=(i_k)_{k=0}^m$ of a fundamental couple the numbering of elements is that of Corollary \ref{2C1} reversed, see \cite[Corollary 3.21 a)]{mu}. Hence in the path associated to $I$ the ES-turns are numbered from right to left. One could avoid this inversion by considering paths from $(\beta,0)$ to $(\alpha,0)$ or by using a different orientation of the diagram, but we prefer our version for typographical reasons.)
\medskip

Conversely, for every lattice path from $(0,\alpha)$ to $(\beta,0)$ not crossing the diagonal the points of the ES-turns can be identified with the gaps in an $\langle \alpha, \beta \rangle$-lean set: 

\begin{lemma}\label{3L3}
Let $\alpha, \beta$ be coprime positive integers. Then there is a bijection between the set of $\langle \alpha,\beta \rangle$-lean sets and the set of lattice paths from $(0,\alpha)$ to $(\beta,0)$ not crossing the diagonal.
\end{lemma}

Therefore counting of $\langle \alpha, \beta \rangle$-lean sets is equivalent to counting of such lattice paths. The latter was considered by Bizley in \cite{biz}. The main idea used there even allows to count the paths with a certain number of ES-turns.
\medskip

The number of all lattice paths with $r$ ES-turns from $(0,\alpha)$ to $(\beta,0)$ is easily computed: The $r$ turning points have $x$-coordinates in the range $\{1, \ldots , \beta -1\}$ and also $y$-coordinates in the range $\{1, \ldots , \alpha -1\}$. Since the sequence of coordinates has to be increasing resp. decreasing there are ${\beta -1\choose r} {\alpha -1\choose r}$ lattice paths. We have to determine how many of these paths stay below the diagonal. To this end we use the concept of a cyclic permutation of a path. 
\medskip

A lattice path with $r$ ES-turns can also be described by a $2 \times (r+1)$-matrix where the $i$-th column contains the numbers of steps downwards and to the right the path takes between the $(i-1)$-th and the $i$-th turning points (where the $0$-th and $(r+1)$-th points are to be understood as $(0,\alpha)$ resp.~$(\beta,0)$). For the path in the example above we get the matrix $
\left(
\begin{array}{cccc}
2 & 1 & 1 & 1\\
1 & 2 & 1 & 3 
\end{array}
\right).
$
\medskip

A \emph{cyclic permutation} of the path is a path belonging to the matrix with cyclically permuted columns. The permuted matrix 
$
\left(
\begin{array}{cccc}
1 & 1 & 2 & 1\\
1 & 3 & 1 & 2 
\end{array}
\right)
$
yields the path
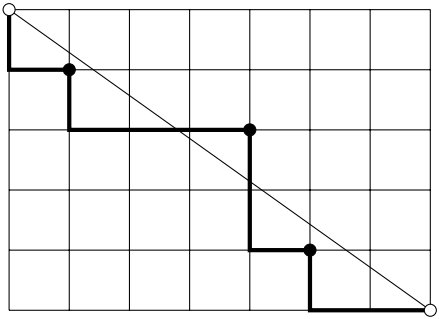
\begin{center}
  \begin{tikzpicture}[scale=0.8]
      \draw[] (0,0) grid [step=1cm](7,5);
      \draw[] (0,5) -- (7,0);
      \draw[ultra thick] (0,5) -- (0,4) -- (1,4) -- (1,3) -- (4,3) -- (4,1) -- (5,1) -- (5,0) -- (7,0);
      \draw[fill=white] (0,5) circle [radius=0.1]; 
      \draw[fill] (1,4) circle [radius=0.1]; 
      \draw[fill] (4,3) circle [radius=0.1]; 
      \draw[fill] (5,1) circle [radius=0.1]; 
      \draw[fill=white] (7,0) circle [radius=0.1]; 
   \end{tikzpicture}
\end{center}

One can also imagine a cyclic permutation in a different way. We extend the path with turning points $P_i=(x_i,y_i)$ beyond $(\beta,0)$ with points $Q_i=(x_i+\beta, y_i-\alpha)$, thus amending a second copy of the original path. The cyclic permutations are the paths from $P_i$ to $Q_i$ with turning points $P_{i+1}, \ldots , P_r, (\beta,0), Q_1, \ldots , Q_{i-1}$:
\begin{center}
  \begin{tikzpicture}[scale=0.63]
     \draw[] (0,0) grid [step=1cm](14,10);
     \draw[ultra thick] (0,10) -- (0,8) -- (1,8) -- (1,7) -- (3,7) -- (3,6) -- (4,6) -- (4,5) -- (7,5);
     \draw[fill=white] (0,10) circle [radius=0.1]; 
     \draw[fill] (1,8) circle [radius=0.1]; 
     \draw[fill=white] (3,7) circle [radius=0.1]; 
     \node [left, above] at (3.4,6.9) {$\scriptstyle P_i$};
     \draw[fill] (4,6) circle [radius=0.1]; 
     \draw[fill] (7,5) circle [radius=0.1]; 

     \draw[ultra thick] (7,5) -- (7,3) -- (8,3) -- (8,2) -- (10,2) -- (10,1) -- (11,1) -- (11,0) -- (14,0);
      \draw[fill] (8,3) circle [radius=0.1]; 
       \draw[fill=white] (10,2) circle [radius=0.1]; 
       \node [left, above] at (10.5,1.9) {$\scriptstyle Q_i$};
       \draw[fill] (11,1) circle [radius=0.1]; 
       \draw[fill=white] (14,0) circle [radius=0.1]; 
       \draw[white] (3,7) -- (10,7) -- (10,2);
       \draw[white] (8,2) -- (3,2) -- (3,6);
      \draw[dashed,thick] (3,7) -- (10,7) -- (10,2)-- (3,2) -- (3,7);
   \end{tikzpicture}
\end{center}
\medskip

A lattice path with $r$ turning points admits $r+1$ cyclic permutations. As we will show now there is always exactly one permutation staying below the diagonal: Again we consider the doubled path as described above. Let $g_0$ denote the line through $P_0:=(0,\alpha)$ and $Q_0:=(\beta,0)$ and $g_i$ the line through $P_i$ and $Q_i$ for $i=1, \ldots, r$. Note that all these lines are parallel and that, since $\alpha$ and $\beta$ are coprime, there is no turning point between $P_i$ and $Q_i$ lying on $g_i$. Hence there is one line $g_j$, $j \in \{0, \ldots , r\}$ namely the one with the greatest distance from the origin, such that the path stays below this line. 

\medskip
\begin{center}
  \begin{tikzpicture}[scale=0.645]
      \draw[ultra thick] (0,10) -- (0,9) -- (3,9) -- (3,7) -- (5,7) -- (5,5) -- (7,5);
      \draw[fill] (0,10) circle [radius=0.1]; 
      \draw[fill] (3,9) circle [radius=0.1]; 
      \draw[fill] (5,7) circle [radius=0.1]; 
       \draw[fill] (7,5) circle [radius=0.1]; 
       \node[above] at (2,7.7) {$ g_0$};
       \node[above] at (8.6,5) {$g_j$};

       \draw[ultra thick] (7,5) -- (7,4) -- (10,4) -- (10,2) -- (12,2) -- (12,0) -- (14,0);
        \draw[fill] (10,4) circle [radius=0.1]; 
        \draw[fill] (12,2) circle [radius=0.1]; 
        \draw[fill] (14,0) circle [radius=0.1]; 
       
         \draw[dashed] (0,10) -- (14,10) -- (14,0) -- (0,0) -- (0,10);
         
          \draw[] (0,10) -- (7,10) -- (7,5) -- (0,5) -- (0,10);
          \draw[] (3,9) -- (10,9) -- (10,4) -- (3,4) -- (3,9);
          \draw[] (5,7) -- (12,7) -- (12,2) -- (5,2) -- (5,7);
          
          \draw[] (0,10) -- (7,5);
          \draw[] (3,9) -- (10,4);
          \draw[] (5,7) --  (12,2);      
   \end{tikzpicture}
\end{center}

Since the line $g_i$ yields the diagonal of the $i$-th cyclic permutation of the original path, indeed exactly one of the permuted paths stays below the diagonal; hence we have proven:

\begin{proposition} 
Let $\alpha$ and $\beta$ be two coprime positive integers. 
\begin{enumerate}
\item[1.] For every lattice path from $(0,\alpha)$ to $(\beta,0)$ there is exactly one cyclic permutation staying below the diagonal.
\item[2.] The number of $\langle \alpha, \beta \rangle$-lean sets with $r$ gaps equals the number of lattice paths with 
$r$ ES-turns from $(0,\alpha)$ to $(\beta,0)$ staying below the diagonal, and this number is given by
\[
\frac{1}{r+1} {\alpha -1 \choose r} {\beta -1 \choose r}.
\]
\end{enumerate}
\end{proposition}

In combination with Corollary \ref{2C2} and Lemma \ref{3L3} this result yields the following theorem:

\begin{theorem}\label{3T1}
Let $\alpha, \beta, r \in \NN$ with $\mathrm{gcd}(\alpha, \beta)=1$. Then the following numbers 
\begin{enumerate}
\item The number of isomorphism classes of $\langle \alpha, \beta \rangle$-semimodules minimally generated by $r+1$ elements.
\item The number of $\langle \alpha, \beta \rangle$-lean sets with $r$ gaps.
\item The number of lattice paths with $r$ ES-turns from $(0,\alpha)$ to $(\beta,0)$ staying below the diagonal.
\end{enumerate}
equal
\[
L_{\alpha,\beta}(r):=\frac{1}{r+1} {\alpha -1 \choose r} {\beta -1 \choose r}.
\]
\end{theorem}

Using standard techniques one can also deduce a formula for $\sum_{r \geq 0} L_{\alpha,\beta}(r)$, recovering results of Bizley, resp.~Beauville, Fantechi--G\"ottsche--van~Straten,  and Piont\-kowski:

\begin{eqnarray*}
\sum_{r \geq 0} L_{\alpha,\beta}(r)& = & \sum_{r \geq 0} \frac{1}{r+1} {\alpha -1\choose r} {\beta -1\choose r} \\
&=& \sum_{r \geq 0} \frac{1}{r+1} {\alpha -1\choose r} \frac{r+1}{\beta} {\beta \choose r+1}\\
&=& \ldots
\end{eqnarray*}
\begin{eqnarray*}
\ldots & = & \frac{1}{\beta} \sum_{r \geq 0} {\alpha -1\choose r} {\beta \choose r+1}\\
& = & \frac{1}{\beta} \sum_{r \geq 1} {\alpha -1\choose r-1} {\beta \choose r} = \frac{1}{\beta} \sum_{r \geq 0} {\alpha -1\choose r-1} {\beta \choose r}\\
&=& \frac{1}{\beta} \sum_{r \geq 0} {\alpha -1\choose \alpha -r} {\beta \choose r}.
\end{eqnarray*}

The Vandermonde convolution yields

\begin{eqnarray*}
\frac{1}{\beta} \sum_{r \geq 0} {\alpha -1\choose \alpha -r} {\beta \choose r} &=&
 \frac{1}{\beta} {\alpha + \beta -1 \choose \alpha}= \frac{1}{\beta} \frac{(\alpha +\beta -1)!}{\alpha ! \cdot (\beta-1)!}= \frac{1}{\alpha+\beta} \frac{(\alpha +\beta)!}{\alpha! \cdot \beta!}\\
&=& \frac{1}{\alpha +\beta} {\alpha +\beta \choose \alpha},
\end{eqnarray*}

which implies the following result.

\begin{theorem}\label{3T2}
Let $\alpha, \beta \in \NN$ be coprime. Then the following numbers
\begin{enumerate}
\item The number of isomorphism classes of $\langle \alpha, \beta \rangle$-semimodules (cf. \cite{bea}, \cite{FGvS}, \cite{pio}).
\item The number of $\langle \alpha, \beta \rangle$-lean sets.
\item The number of lattice paths from $(0,\alpha)$ to $(\beta,0)$ staying below the diagonal (cf. \cite{biz}).
\end{enumerate}
equal
\[
L_{\alpha,\beta}:=\sum_{r \geq 0} L_{\alpha,\beta}(r) =\frac{1}{\alpha +\beta} {\alpha +\beta \choose \alpha}.
\]
\end{theorem}

\begin{remark}
In particular $L_{\alpha,\beta}$ gives the number of $(\alpha,\beta)$-fundamental couples and hence the number of inequalities appearing in condition $(\star)$ of Theorem 3.13 in \cite{mu}, also see \cite[Remark 3.12]{mu}.
\end{remark}

\begin{remark}
In the special case $\beta=\alpha+1$ the numbers $L_{\alpha,\beta}(r)$ and $L_{\alpha,\beta}$ coincide with certain combinatorial numbers. We get
\[
L_{\alpha,\alpha+1}(r)=\frac{1}{\alpha}{\alpha \choose r}{\alpha \choose r+1} = N(\alpha,r+1),
\]
the so-called \emph{Narayana number} (cf. \cite{pro}). Moreover, $L_{\alpha,\alpha+1}$ agrees with the \emph{Catalan number} $C_{\alpha}$ since
\[
L_{\alpha,\alpha+1}=\sum_{r\geq 0} L_{\alpha, \alpha+1}(r)=\sum_{r \geq 0} N(\alpha, r+1)= C_{\alpha},
\]
see \cite[Exercise 6.36]{stan}.
\end{remark}

\section{Syzygies of $\langle \alpha, \beta \rangle$-semimodules and lattice paths} \label{section4}

We consider now the sequences appearing in the second position of a fundamental couple. Let $[I,J]$ be a fundamental couple with sequences $I=[i_0=0, \ldots , i_n]$ and $J=[j_0,\ldots,j_n]$. By definition, the elements $j_1,\ldots ,j_{n-1}$ are gaps of $\langle \alpha, \beta \rangle$ such that
\[
j_k \equiv i_k ~ \mbox{~mod~} \alpha ~~ \mathrm{~and~} j_k \equiv i_{k+1} ~ \mbox{~mod~} \beta.
\]
An inspection of the lattice path belonging to $I$ shows that these gaps $j_1, \ldots , j_{n-1}$ correspond to the inner SE-turning points of the path. By extension of the labeling beyond the axis we can even identify $j_0$ and $j_n$ with the remaining SE-turns. For illustration see again the example of the previous section:
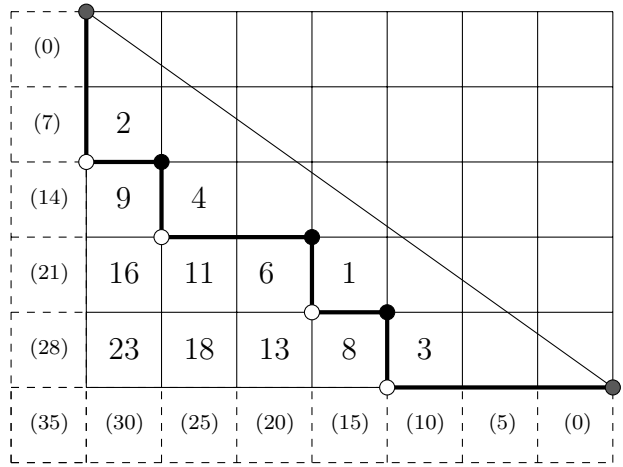
\vspace{12pt}

\begin{center}
  \begin{tikzpicture}
    \draw[] (0,0) grid [step=1cm](7,5);
    \draw[] (0,5) -- (7,0);
    \draw[ultra thick] (0,5) -- (0,3) -- (1,3) -- (1,2) -- (3,2) -- (3,1) -- (4,1) -- (4,0) -- (7,0);
     \node [below right] at (0.15,0.8) {$23$};
      \node [below right] at (1.15,0.8) {$18$};
       \node [below right] at (2.15,0.8) {$13$};
       \node [below right] at (3.25,0.8) {$8$};
       \node [below right] at (4.25,0.8) {$3$};
       
       \node [below right] at (0.15,1.8) {$16$};
      \node [below right] at (1.15,1.8) {$11$};
       \node [below right] at (2.15,1.8) {$6$};
       \node [below right] at (3.25,1.8) {$1$};
       
           \node [below right] at (0.25,2.8) {$9$};
      \node [below right] at (1.25,2.8) {$4$};
      
          \node [below right] at (0.25,3.8) {$2$};
          
         \draw[dashed] (-1,-1) grid [step=1cm](7,0);
           \draw[dashed] (-1,-1) grid [step=1cm](0,5);
  
      \node [below right] at (-0.9,-0.2) {$\scriptstyle (35)$};
      \node [below right] at (0.1,-0.2) {$\scriptstyle (30)$};
       \node [below right] at (1.1,-0.2) {$\scriptstyle (25)$};
        \node [below right] at (2.1,-0.2) {$\scriptstyle (20)$};
         \node [below right] at (3.1,-0.2) {$\scriptstyle (15)$};
          \node [below right] at (4.1,-0.2) {$\scriptstyle (10)$};
           \node [below right] at (5.2,-0.2) {$\scriptstyle (5)$};
            \node [below right] at (6.2,-0.2) {$\scriptstyle (0)$};

        \node [below right] at (-0.9,0.8) {$\scriptstyle (28)$};
         \node [below right] at (-0.9,1.8) {$\scriptstyle (21)$};
         \node [below right] at (-0.9,2.8) {$\scriptstyle (14)$};
         \node [below right] at (-0.85,3.8) {$\scriptstyle (7)$};
        \node [below right] at (-0.85,4.8) {$\scriptstyle (0)$};

    \draw[fill=black!65] (0,5) circle [radius=0.1]; 
    \draw[fill] (1,3) circle [radius=0.1]; 
     \draw[fill] (3,2) circle [radius=0.1]; 
     \draw[fill] (4,1) circle [radius=0.1]; 
      \draw[fill=black!65] (7,0) circle [radius=0.1];

        \draw[fill=white] (0,3) circle [radius=0.1]; 
     \draw[fill=white] (1,2) circle [radius=0.1]; 
     \draw[fill=white] (3,1) circle [radius=0.1]; 
      \draw[fill=white] (4,0) circle [radius=0.1]; 
   \end{tikzpicture}
\end{center}
\begin{center}
{\small $I=[0,8,6,9]$ and $J=[15,13,16,14]$.}
\end{center}
\medskip
\vspace{12pt}

Next we explain the meaning of $J$ in terms of $\langle \alpha, \beta \rangle$-semimodules: Every $\langle \alpha, \beta \rangle$-semimodule $\Delta$ yields another $\langle \alpha, \beta \rangle$-semimodule $\mathrm{Syz}(\Delta)$.

\begin{definition}\label{4D1}
Let $I$ be an $\langle \alpha, \beta \rangle$-lean set, and
let $\Delta$ be the $\langle \alpha, \beta \rangle$-semimodule generated by $I$. The \emph{syzygy} of $\Delta$ is the $\langle \alpha, \beta \rangle$-semimodule 
\[
\mathrm{Syz}(\Delta):=\bigcup_{\substack{i,i' \in I\\i \neq i'}} \left ( \big(i+\langle \alpha, \beta \rangle \big) \cap \big(i' + \langle \alpha, \beta \rangle \big)  \right ).
\]
\end{definition}

The semimodule $\mathrm{Syz}(\Delta)$ consists of those elements in $\Delta$ which admit more than one presentation of the form $i + x$ with $i \in I, x \in \langle \alpha, \beta \rangle$. The name \emph{syzygy} may be justified by considering an analogue of $\Delta$ in the setting of commutative algebra: 
\newpage
Let $R=\FF[t^{\alpha},t^{\beta}]$, then $\Delta$ can be identified with the $R$-submodule $M$ of $\FF[t]$ generated by $\{t^{i} \mid i\in I \}$. Let
\begin{eqnarray*}
\bigoplus_{i \in I} R(-i) &\stackrel{\varphi}{\longrightarrow}& M \\
(f_i) & \mapsto & \sum_i f_i t^{i}
\end{eqnarray*}
be the first step in a graded minimal free resolution of $M$. By \cite[Lemma 2.3]{pio} the kernel of $\varphi$ is generated by (homogeneous) elements $v_{i,i'}$ of the form 
$$
(0, \ldots , t^{\gamma_i}, 0, \ldots , 0, -t^{\gamma_{i'}},0, \ldots , 0),
$$
with the non-zero entries in positions $i,i' \in I$. Since 
$$
\deg(v_{i,i'})=\gamma_i+i = \gamma_{i'}+i' \in \left ( i+\langle \alpha, \beta \rangle \right ) \cap \left ( i'+\langle \alpha, \beta \rangle \right )
$$
the module $\ker \varphi$ ist non-zero exactly in the degrees   contained in $\mathrm{Syz}(\Delta)$. 
\medskip

The connection between fundamental couples and syzygies is described in the following theorem:

\begin{theorem}\label{4T2}
Let $[I,J]$ be an $\langle \alpha, \beta \rangle$-fundamental couple and let $\Delta$ be the $\langle \alpha, \beta \rangle$-semimodule generated by the elements of $I$. Then 
\[
\mathrm{Syz}(\Delta)=\bigcup_{0 \leq k < m \leq n} \Big( \left ( i_k + \langle \alpha, \beta \rangle \right ) \cap \left (  i_{m}+\langle \alpha, \beta \rangle \right )  \Big )=\bigcup_{k=0}^{n}(j_k + \langle \alpha, \beta \rangle).
\]
\end{theorem}

\proof

By definition of a fundamental couple we have $j_k=i_k+r\alpha=i_{k+1}+s\beta$ with some $r,s \in \NN$, hence the inclusion $\supseteq$ is clear.
In order to show the other inclusion, we consider $i_k + \langle \alpha, \beta \rangle \cap i_{m}+\langle \alpha, \beta \rangle$. For every $\gamma_m \in \langle \alpha, \beta \rangle$ we have 
$i_m+\gamma_m \in i_k + \langle \alpha, \beta \rangle$ if and only if $i_m-i_k+\gamma_m \in \langle \alpha, \beta \rangle$. As mentioned in the previous section, $i_k$ and $i_m$ can be written in the form
\[
i_k=\alpha \beta -a_k \alpha -b_k \beta, ~~~i_m=\alpha \beta -a_m \alpha -b_m \beta;
\]
by Lemma \ref{3L2} we may assume $a_k >a_m$ and $b_k<b_m$. Since
\[
i_m-i_k = (a_k -a_m)\alpha +(b_k-b_m)\alpha=\alpha \beta -(\beta-a_k+a_m)\alpha -(b_m-b_k)\beta
\]
the characterization of $\ZZ \setminus \langle \alpha, \beta \rangle$ in Lemma \ref{3L1} implies
\[
\{\gamma \in \langle \alpha, \beta \rangle \mid \gamma +i_m-i_k \in \langle \alpha, \beta \rangle\}=\left ( (\beta -a_k+a_m)\alpha +\langle \alpha, \beta \rangle \right ) \cup \left (  (b_m-b_k)\beta +\langle \alpha, \beta \rangle \right ).
\]
This means
\begin{eqnarray*}
&&\left ( i_m + \langle \alpha, \beta \rangle \right ) \cap \left (  i_k+\langle \alpha, \beta \rangle \right ) \\
&=& \left ( i_m+(\beta-a_k+a_m)\alpha + \langle \alpha, \beta \rangle \right ) \cup \left (  i_m + (b_m-b_k)\beta + \langle \alpha, \beta \rangle \right ).
\end{eqnarray*}

Moreover, 
\begin{eqnarray*}
i_m + (\beta -a_k+a_m)\alpha &=& \alpha\beta -a_k \alpha - b_m \beta + \alpha \beta\\
&=& (\beta-a_1)\alpha + (a_1-a_k)\alpha + (\alpha-b_m)\beta\\
&=& j_0 + (a_1-a_k)\alpha + (\alpha-b_m)\beta \in j_0+ \langle \alpha, \beta \rangle,
\end{eqnarray*}
and on the other hand
\begin{eqnarray*}
i_m + (b_m-b_k)\beta &=& \alpha\beta - a_m\alpha - b_m\beta + (b_m-b_k)\beta\\
&=&\alpha \beta - a_m \alpha - b_k \beta\\
&=& \alpha \beta -a_{k+1} \alpha - b_k \beta + (a_{k+1}-a_m)\alpha\\
&=& j_k + (a_{k+1}-a_m)\alpha \in j_k + \langle \alpha, \beta \rangle,
\end{eqnarray*}
hence $i_m + \langle \alpha, \beta \rangle \cap i_k + \langle \alpha, \beta \rangle \subseteq j_0 + \langle \alpha, \beta \rangle \cup j_k + \langle \alpha, \beta \rangle$. \qed

\begin{corollary}\label{4C1}
Let $[I,J]$ be an $\langle \alpha, \beta \rangle$-fundamental couple and let $\Delta $ be the $\langle \alpha, \beta \rangle$-semimodule generated by the elements of $I$. We have 
\[
\mathrm{Syz}(\Delta)=\bigcup_{k=0}^{n-1} \Big( \left ( i_k + \langle \alpha, \beta \rangle \right ) \cap \left ( i_{k+1}+\langle \alpha, \beta \rangle \right ) \Big ) \cup \big (\left ( i_0+\langle \alpha, \beta \rangle \right ) \cap \left ( i_n +\langle \alpha, \beta \rangle \right ) \big ).
\]
\end{corollary}

\proof
This follows immediately from $j_k \in i_k + \langle \alpha, \beta \rangle \cap i_{k+1}+\langle \alpha, \beta \rangle$ for $k=1, \dots , n-1$ resp.~$j_n \in i_0+\langle \alpha, \beta \rangle \cap i_n + \langle \alpha, \beta \rangle$ and the previous theorem.
\qed

\section{Orbits}

Let $[I,J]$ a fundamental couple and let
\[
\left(
\begin{array}{cccc}
y_0 & y_1 & \ldots & y_n   \\
x_0 & x_1 & \ldots & x_n 
\end{array}
\right)
\]
be the matrix describing the path for the semimodule $\Delta$ generated by $I$.
We consider a second lattice path from $(0,b_n)$ (the point associated to $j_n$) to $(\beta, b_n-\alpha)$ with ES-turns in the SE-turning points of the first path (those points representing $j_{n-1}, \ldots , j_{0}$). The matrix for this path is given by
\[
\left(
\begin{array}{ccccc}
y_1 & y_2 & \ldots & y_n& y_0   \\
x_ 0& x_1 & \ldots & x_{n-1} & x_n 
\end{array}
\right).
\]
It is easily seen that---up to a cyclic permutation of columns---this matrix also describes the path belonging to the normalization $\Delta^{\circ}$ of $\Delta$: In terms of lattice paths normalizing $\Delta$ means translation of the path such that the ES-turn belonging to $\min J$ is moved to $(0,\alpha)$ and the part of the path left of this point will be appended behind the former end point. 
\medskip

The procedure of building a syzygy can be iterated; we set 
\[
\mathrm{Syz}^{(k)}(\Delta):=\mathrm{Syz}(\mathrm{Syz}^{(k-1)}(\Delta)), ~~ k \geq 2.
\]
From the matrix description of the path for $\mathrm{Syz}(\Delta)$ it is clear that $\mathrm{Syz}^{(n+1)}(\Delta)\cong \Delta$. Now we consider under which conditions even lower syzygies of $\Delta$ are isomorphic with $\Delta$.

\begin{definition} 
A sequence $\Delta_1, \ldots , \Delta_{\ell}$ of distinct  $\langle \alpha, \beta \rangle$-semimodules such that $\Delta_k=\mathrm{Syz}(\Delta_{k-1})^{\circ}$ for $k=2,\ldots, \ell$, and $\Delta_1=\mathrm{Syz}(\Delta_{\ell})^{\circ}$ is called an orbit of length $\ell$ or, for short, an $\ell$-orbit (of $\langle \alpha, \beta \rangle$-semimodules). The $1$-orbits, i.~e.~$\langle \alpha, \beta \rangle$-semimodules $\Delta$ with $\Delta \cong \mathrm{Syz}(\Delta)$, are called $\langle \alpha, \beta \rangle$-fixed points.
\end{definition}

We want to count the number of $\ell$-orbits of $\langle \alpha, \beta \rangle$-semimodules with $n$ generators. To this end we investigate the structure of a semimodule $\Delta$ with $\mathrm{Syz}^{(\ell)}(\Delta)\cong \Delta$; we may restrict our attention to the case of $\ell$ dividing $n$, since the length of an orbit can be viewed as the order of an element in the cyclic group of order $n$.  Let
\[
\left(
\begin{array}{cccc}
y_0 & y_1 & \ldots  & y_{n-1} \\
x_0 & x_1 & \ldots & x_{n-1} 
\end{array}
\right)
\]
be the matrix for the path belonging to $\Delta$. Then, as mentioned above, this matrix and the matrix
\[
\left(
\begin{array}{ccccc}
y_{\ell} & y_{\ell +1} & \ldots &y_{\ell-2}&y_{\ell-1}  \\
x_0 & x_1 & \ldots &x_{n-1}& x_n
\end{array}
\right),
\]
have to be equal up to cyclic permutation of columns. 
This means that there exists a $k \in \{1, \ldots , n-1\}$ such that
\[
\left(
\begin{array}{ccc}
y_{k+\ell} & y_{k+\ell + 1} & \ldots \\
x_k & x_{k+1} & \ldots   
\end{array}
\right)=
\left(
\begin{array}{ccc}
y_0 & y_1 & \ldots   \\
x_0 & x_1 & \ldots  
\end{array}
\right),
\]
we may assume that $k$ is minimal with this property. Since 
\[
\left(
\begin{array}{ccc}
y_{k+2\ell} & y_{k+2\ell +1} & \ldots \\
x_k & x_{k+1} & \ldots  
\end{array}
\right)=
\left(
\begin{array}{ccc}
y_{\ell} & y_{\ell + 1} & \ldots   \\
x_0 & x_1 & \ldots
\end{array}
\right)
\]
we have 
\[
\left(
\begin{array}{ccc}
y_{2k+2\ell} & y_{2k+2\ell +1} & \ldots   \\
x_{2k} & x_{2k+1} & \ldots  
\end{array}
\right)=
\left(
\begin{array}{ccc}
y_{k+\ell} & y_{k+\ell + 1} & \ldots \\
x_k & x_{k+1} & \ldots
\end{array}
\right)=
\left(
\begin{array}{ccc}
y_0 & y_1 & \ldots   \\
x_0 & x_1 & \ldots 
\end{array}
\right).
\]
By induction we get $x_j=x_{rk+j}$ for $j \in \{0, \ldots k-1\}$, $r \geq 0$. Hence the bottom row of the matrix is of the form
\[
[x_0~ \ldots ~x_{k-1}] ~ [x_0 ~\ldots ~x_{k-1}] ~ \ldots  ~ [x_0 ~\ldots ~x_{k-1}].
\]
On the other hand we can find a minimal $m \in \{2, \ldots, n\}$ such that 
\[
\left(
\begin{array}{ccc}
y_{m-\ell} & y_{m-\ell+1} & \ldots   \\
x_{m} & x_{m+1} & \ldots  
\end{array}
\right)=
\left(
\begin{array}{ccc}
y_0 & y_1 & \ldots   \\
x_0 & x_1 & \ldots  
\end{array}
\right),
\]
and with the same reasoning as above, the top row of the matrix is of the form
\[
[y_0~ \ldots ~y_{m-1}] ~ [y_0 ~\ldots ~y_{m-1}] ~ \ldots  ~ [y_0 ~\ldots ~y_{m-1}].
\]
Therefore the matrix for $\Delta$ looks like
\[
\left(
\begin{array}{lcr}
\mbox{[}y_0~ \ldots  \ldots ~y_{m-1}] &  \ldots  & [y_0 ~\ldots  \ldots~y_{m-1}] \\
\mbox{[}x_0~ \ldots ~x_{k-1}] & \ldots & [x_0 ~\ldots ~x_{k-1}]  
\end{array}
\right)
\]
with $m'$ blocks $[y_0 ~\ldots ~y_{m-1}]$ and $k'$ blocks $[x_0~ \ldots ~x_{k-1}] $. Since 
\begin{equation} \label{5E1}
m' \cdot \sum_{j=0}^{m-1} y_j=\alpha ~ ~ ~ ~\mbox{~~~and~~~} ~ ~ ~ ~k'\cdot \sum_{j=0}^{k-1}x_j=\beta,
\end{equation}
the numbers $m'$ and $k'$ divide $\alpha$ resp.~$\beta$, so in particular they are coprime. By assumption, $\ell$ is the least positive integer $p$ such that the matrices 
\[
\left(
\begin{array}{ccc}
y_0 & y_1 & \ldots   \\
x_0 & x_1 & \ldots  
\end{array}
\right)
~~\mbox{~~and~~}~~
\left(
\begin{array}{ccc}
y_{p} & y_{p + 1} & \ldots   \\
x_0 & x_1 & \ldots  
\end{array}
\right)
\]
contain the same columns. This is the case if and only if there are $r,s \in \mathbb{N}$ with $p+s m=r k$. This implies that $p$ has to be contained in the ideal generated by $k$ and $m$. Hence, by minimality of $\ell$, we get $\ell=\mathrm{gcd}(k,m)$, and so we may write $k=\tilde{k}\ell$ and $m=\tilde{m}\ell$. From $kk'=n=mm'$ we get $\tilde{k}k'=\tilde{m}m'=\frac{n}{\ell}$; by $\mathrm{gcd}(\tilde{k},\tilde{m})=1$ this implies $\tilde{k}=m'$ and $\tilde{m}=k'$, and by $k' \mid \alpha$ and $m' \mid \beta$ moreover $k=\mathrm{gcd}(\beta,\frac{n}{\ell})$ and $m=\mathrm{gcd}(\alpha,\frac{n}{\ell})$, in particular $\frac{n}{\ell} \mid \alpha \beta$. By now we have shown the following proposition:

\begin{proposition}\label{5L1}
Let $\Delta$ be an $\langle \alpha, \beta \rangle$-semimodule with $n$ generators. If $\Delta$ is an element of an $\ell$-orbit, then $\ell\!\mid \!n$ and $\frac{n}{\ell} \!\mid\!\alpha \beta$. The corresponding matrix is of the form
\begin{equation}\label{5E2}
\left(
\begin{array}{lcr}
[y_0~ \ldots ~y_{m-1}] & \ldots & [y_0 ~\ldots ~y_{m-1}]  \\
\mbox{$[$}x_0~ \ldots  \ldots ~x_{k-1}] &  \ldots  & [x_0 ~\ldots  \ldots~x_{k-1}] 
\end{array}
\right),
\end{equation}
where $k=\ell \cdot \mathrm{gcd}(\alpha,\frac{n}{\ell})$ and $m=\ell \cdot \mathrm{gcd}(\beta,\frac{n}{\ell})$.
\end{proposition}

In fact the structure described in the previous proposition is shared by all  $\langle \alpha, \beta \rangle$-semimodules $\Delta$ with $n$ generators and $\mathrm{Syz}^{(\ell)}(\Delta)\cong \Delta$---not only by elements of $\ell$-orbits: 

\begin{lemma}\label{5L2}
Let $\Delta$ be an $\langle \alpha, \beta \rangle$-semimodule with $n$ generators, let $\ell \in \mathbb{N}$ be a divisor of $n$. If $\mathrm{Syz}^{(\ell)}(\Delta)\cong \Delta$ then the matrix for $\Delta$ can be written in the form mentioned in Proposition \ref{5L1}. 
\end{lemma}

\proof 
Let $d=\min \{t \in \mathbb{N} \mid \mathrm{Syz}^{(t)}(\Delta)\cong \Delta \}$. Then $d \mid \ell$, and we only have to consider the case $d < \ell$. By Proposition \ref{5L1} the matrix for $\Delta$ can be written in the form 
$$
\left(
\begin{array}{lcr}
[y_0~ \ldots ~y_{m_d-1}] & \ldots & [y_0 ~\ldots ~y_{m_d-1}]  \\
\mbox{$[$}x_0~ \ldots  \ldots ~x_{k_d-1}] &  \ldots  & [x_0 ~\ldots  \ldots~x_{k_d-1}] 
\end{array}
\right)
$$
\medskip
with $k'_d=\mathrm{gcd}(\beta, \frac{n}{d})$ blocks $[x_0~  \ldots ~x_{k_d-1}]$ and  $m'_d=\mathrm{gcd}(\alpha, \frac{n}{d})$ blocks $[y_0~  \ldots ~y_{m_d-1}]$. We want to obtain a matrix with $k'_{\ell}=\mathrm{gcd}(\beta, \frac{n}{\ell})$ blocks $[x_0~  \ldots ~x_{k_{\ell}-1}]$ and  $m'_{\ell}=\mathrm{gcd}(\alpha, \frac{n}{\ell})$ blocks $[y_0~  \ldots ~y_{m_{\ell}-1}]$. Since $k'_{\ell} \mid k'_{d}$ and $m'_{\ell} \mid m'_{d}$, this is easily done by concatenating $\frac{k'_{d}}{k'_{\ell}}$ $x$-blocks resp. $\frac{m'_{d}}{m'_{\ell}}$ $y$-blocks. 
\qed

\begin{corollary}\label{5C0}
Let $\ell \in \mathbb{N}$ be a divisor of $n <\alpha$. Then there exists an $\ell$-orbit of $\langle \alpha, \beta \rangle$-semimodules with $n$ generators.
\end{corollary}

\proof 
Let $k,k',m,m'$ be as in the proof of Prop. \ref{5L1}. Then $\sum y_i = \frac{\alpha}{m'}=m\frac{\alpha}{n} >m$ resp. $\sum x_i > k$. Hence $r:= \frac{\alpha}{m'}-m+1>1$ and $s:= \frac{\beta}{k'}-k+1>1$, so we may choose the matrix (\ref{5E2}) to be built of blocks $[y_0 \ldots]=[1~1\ldots 1~r]$, $[x_0 \ldots]=[1~ 1\ldots 1~s]$. Since these blocks cannot be split into smaller ones, the corresponding semimodule cannot be part of a $d$-orbit with $d <\ell$.
\qed
\medskip

Lemma \ref{5L2} allows to count $\langle \alpha, \beta \rangle$-semimodules $\Delta$ with $\mathrm{Syz}^{(\ell)}(\Delta)\cong \Delta$ and $n$ generators. Using the same notation as in the deduction of Proposition \ref{5L1}, we find the following:  By (\ref{5E1}) the steps in the $y$-block sum up to $\frac{\alpha}{m'}$, therefore the partial sums $\sum_{j=0}^{r} y_j$ for $r=0, \ldots , m-2$ (these are the $y$-coordinates of the first $m-1$ ES-turns in the corresponding path) have to be chosen in the range $1, \ldots , \frac{\alpha}{m'}-1$. Hence there are ${\frac{\alpha}{m'}-1 \choose m-1}$ different $y$-blocks. Similarly there are ${\frac{\beta}{k'}-1 \choose k-1}$ different $x$-blocks. Any combination of an $x$- and a $y$-block yields a matrix of the form (\ref{5E2}), and so there are ${\frac{\alpha}{m'}-1 \choose m-1}{\frac{\beta}{k'}-1 \choose k-1}$ of them. But, as in the counting of lattice paths in section 3, only one of the $n$ cyclic permutations of the matrix is admissible, hence we have to divide by $n$: 

\begin{theorem}\label{5T1}
There are
\begin{equation}\label{5E4}
\frac{1}{n}{\frac{\alpha}{\mathrm{gcd}(\frac{n}{\ell},\alpha)}-1\choose \ell \cdot \mathrm{gcd}(\frac{n}{\ell},\beta)-1}{\frac{\beta}{\mathrm{gcd}(\frac{n}{\ell},\beta)}-1\choose \ell \cdot \mathrm{gcd}(\frac{n}{\ell},\alpha)-1}
\end{equation}
$\langle \alpha, \beta \rangle$-semimodules $\Delta$ with $\mathrm{Syz}^{(\ell)}(\Delta)\cong \Delta$ generated by $n$ elements.
\end{theorem}

In particular we get a formula for the number of $\langle \alpha, \beta \rangle$-fixed points:

\begin{corollary}\label{5C1}
For any integer $n \leq \alpha$ with $n \mid \alpha \beta$ there are
\begin{equation}\label{5E3}
\frac{1}{n}{\frac{\alpha}{\mathrm{gcd}(n,\alpha)}-1\choose \mathrm{gcd}(n,\beta)-1}{\frac{\beta}{\mathrm{gcd}(n,\beta)}-1\choose \mathrm{gcd}(n,\alpha)-1}
\end{equation}
$\langle \alpha, \beta \rangle$-fixed points with $n$ generators.
\end{corollary}

\begin{remark}
\noindent 
\begin{enumerate}
\item In the case of $n=\ell$ the number provided by Theorem \ref{5T1} agrees with that of all $\langle \alpha, \beta \rangle$-semimodules with $n$ generators (Thm.~\ref{3T1}), in accordance to the fact that $\mathrm{Syz}^{(n)}(\Delta)\cong \Delta$ for all these semimodules.
\item For $n=\alpha$ the formula (\ref{5E3}) yields $\frac{1}{\alpha}{\beta-1\choose \alpha-1}$, which equals $L_{\alpha,\beta}(\alpha-1)$. Hence \emph{all} $\langle \alpha, \beta \rangle$-semimodules with maximal number of generators are fixed points---in this case there are only entries ``$1$" in the top row of the matrix. Note that in the case of $\alpha$ and $\beta$ being prime numbers there are no other fixed points. 
\end{enumerate}
\end{remark}

One has to keep in mind that formula (\ref{5E4}) does not only count the elements with $n$ generators in the $\ell$-orbits, but \emph{all} semimodules with $n$ generators such that $\mathrm{Syz}^{(\ell)}(\Delta)\cong \Delta$,
including those with $\mathrm{Syz}^{(d)}(\Delta)\cong \Delta$ and $d < \ell$, in particular the fixed points. However, it is possible to compute the number of $\ell$-orbits using the inclusion-exclusion principle, see the next and closing example.

\begin{example}
Let $\langle \alpha,\beta \rangle =\langle 15,16 \rangle$. We want to compute how many orbits consisting of semimodules with $12$ generators exist. Since $12$ divides $15\cdot 16$ there are $\ell$-orbits for each divisor $\ell$ of $12$. Denote the set of elements of those orbits by $\mathrm{Orb}_{\ell}$ and the set of those semimodules $\Delta$ with $n$ generators and $\mathrm{Syz}^{(\ell)}(\Delta)\cong \Delta$ by $A_{\ell}$. By Theorem \ref{5T1} we get
\[
|A_{1}|=1,~~|A_2|=7,~~|A_3|=91,~~|A_4|=455,~~|A_6|=637,~~|A_{12}|=41405.
\]
We have $\mathrm{Orb}_{1}=A_1$, and one easily checks that
\begin{eqnarray*}
\mathrm{Orb}_{2}&=&A_2 \setminus A_1\\
\mathrm{Orb}_{3}&=&A_3 \setminus A_1\\
\mathrm{Orb}_{4}&=&A_4 \setminus A_2\\
\mathrm{Orb}_{6}&=&A_6 \setminus (A_2 \cup A_3)\\
\mathrm{Orb}_{12}&=&A_{12} \setminus \cup_{i=1,2,3,4,6} \mathrm{Orb}_{i}.
\end{eqnarray*}
Since $A_2 \cap A_3=A_1$ we get
\begin{eqnarray*}
|\mathrm{Orb}_{1}|&=&1\\
|\mathrm{Orb}_{2}|&=&|A_2| - |A_1|=6\\
|\mathrm{Orb}_{3}|&=&|A_3| - |A_1|=90\\
|\mathrm{Orb}_{4}|&=&|A_4|-|A_2|=448\\
|\mathrm{Orb}_{6}|&=&|A_6 \setminus (A_2 \cup A_3)|=|A_6|-|A_2|-|A_3|+|A_1|=540\\
|\mathrm{Orb}_{12}|&=&41405-540-448-90-6-1=40320,
\end{eqnarray*}
and thus the following numbers of $\ell$-orbits:
\begin{center}
\begin{tabular}{c||c|c|c|c|c|c}
$\ell$ & 1 & 2 & 3 & 4 & 6 &12\\\hline
$\frac{1}{\ell} |\mathrm{Orb}_{\ell}|$ & 1& 3 & 30 & 112& 90 & 3360
\end{tabular}
\end{center}

Concluding the discussion of this example we determine the single fixed point. By Proposition \ref{5L1} the rows of the corresponding matrix have to consist of three blocks $[y_0, y_1,y_2,y_3]$ resp.~four blocks $[x_0,x_1,x_2]$. Since the entries in those blocks sum up to $\frac{15}{3}=5$ resp.~$\frac{16}{4}=4$ both blocks have to contain entries ``$1$'' and a single entry ``$2$'' i.e.~up to cyclic permutation the matrix of the fixed point looks like
\begin{equation*}
\left(
\begin{array}{cccccccccccc}
1 & 1& 1&2 &1 &1 &1 &2 &1 &1 &1 &2   \\
1 &1 &2 &1 &1 &2 &1 &1 &2 &1 &1 &2 
\end{array}
\right).
\end{equation*}
By considering the corresponding lattice path one can deduce that the admissible permutation of the path is described by
\begin{equation*}
\left(
\begin{array}{cccccccccccc}
2 & 1& 1&1 &2 &1 &1 &1 &2 &1 &1 &1   \\
1 &1 &2 &1 &1 &2 &1 &1 &2 &1 &1 &2 
\end{array}
\right),
\end{equation*}
and the minimal set of generators of the $\langle 15,16\rangle$-semimodule $\Delta$ belonging to this path is given by
$
I=\{0,14,13,12,10,9,8,22,5,4,18,17\};
$
the first syzygy of $\Delta$ is generated by
$
J=\{30,29,28,42,25,24,38,37,20,34,33,32\}.
$
\end{example}

\

\subsection*{Acknowledgements}
The authors wish to express their gratitude to Evgeny Gorsky and Carlos Mariju\'an L\'opez for many helpful suggestions as well as for their encouraging interest in our work.

The first author was partially supported by the Spanish Government Ministerio de Educaci\'on y Ciencia, grant MTM2007-64704 and Ministerio de Econom\'ia y Competitividad, grant MTM2012--36917--C03--03, in cooperation with the European Union in the framework of the founds ``FEDER''.


\begin{thebibliography}{99}

\bibitem{bea} A.~Beauville, \emph{Counting rational curves on K3 surfaces}. Duke Math. J. \textbf{97}~(1999), 99--108.

\bibitem{biz} M.~T.~L.~Bizley, \emph{Derivation of a new formula for the number of minimal lattice paths from $(0,0)$ to $(km,kn)$ having just $t$ contacts with the line $my=nx$ and having no points above this line; and a proof of Grossman's formula for the number of paths which may touch but do not rise above this line}. J. Inst. Actuar. \textbf{80} (1954) 55--62.

\bibitem{FGvS} B.~Fantechi, L.~G\"ottsche, D.~van~Straten, \emph{Euler number of the compactified Jacobian and multiplicity of rational curves}. J. Alg. Geom. \textbf{8} (1999), 115--133.

\bibitem{mu} J.~J.~Moyano-Fern\'andez, J.~Uliczka, \emph{Hilbert depth of graded modules over polynomial rings in two variables}. J.~of~Algebra \textbf{373} (2013)~130--152.

\bibitem{pio} J.~Piontkowski, \emph{Topology of the CompactiÞed Jacobians of Singular Curves}. Math. Z. \textbf{255}(1) (2007), 195--226.

\bibitem{pro} H.~Prodinger, \emph{A correspondence between ordered tress and noncrossing partitions}. Discrete Mathematics \textbf{46} (1983), 205--206.

\bibitem{rosales} J.~C.~Rosales, \emph{Fundamental gaps of numerical semigroups generated by two elements}.
                  Linear Algebra and its Appl.~\textbf{405}~(2005)~200--208.
                  
\bibitem{stan} R.~P.~Stanley, \emph{Enumerative Combinatorics} Vol.~2. Cambridge Studies in Adv. Math. 62, Cambridge U.P. 1999

\end{thebibliography}
\end{document}